\documentclass[10pt]{article}
\usepackage[tbtags]{amsmath}
\usepackage{epsfig,amstext,amssymb,amsthm,latexsym}
\catcode`\@=11
\@addtoreset{equation}{section}
\catcode`\@=12

\usepackage{amssymb,color}
\definecolor{c20}{rgb}{0.,0.7,0.}
\definecolor{c30}{rgb}{0.,0.,1.}
\definecolor{c40}{rgb}{1,0.1,0.7}
\definecolor{c50}{rgb}{1,0,0}

\def\cre#1{#1}
\def\resub#1{ \textcolor{c20}{#1}}
\def\resub#1{ #1}

\allowdisplaybreaks[4]
\newcommand{\nwc}{\newcommand}
\nwc{\COM}[1]{}

\nwc{\vs}[1]{\vskip #1 cm}

\newtheorem{theo}{Theorem}[section]
\newtheorem{sat}[theo]{Proposition}
\newtheorem{de}[theo]{Definition}
\newtheorem{lem}[theo]{Lemma}

\newtheorem{korr}[theo]{Corollary}
\newtheorem{remark}[theo]{Remark}
\newtheorem{exxa}[theo]{Example}

\newcommand{\nelem}[1]{{Lemma \ref{#1}}}

\newcommand{\netheo}[1]{{Theorem \ref{#1}}}

\newcommand{\kb}[1]{\boldsymbol{#1}}
\newcommand{\vk}[1]{\kb{#1}}

\def\FRE{\mbox{Fr\'{e}chet }}
\def\a{\vk{a}}

\def\x{\vk{x}}
\def\y{\vk{y}}
\def\X{\vk{X}}

\newcommand{\ve}{\varepsilon}

\newcommand{\abs}[1]{\lvert #1 \rvert}
\newcommand{\Abs}[1]{ \Bigl \lvert #1 \Bigr \rvert}

\newcommand{\norm}[1]{\lVert #1 \rVert}

\newcommand{\E}[1]{\mbox{\rm$\vk{E}$}\{#1\}}

\newcommand{\pk}[1]{\mbox{\rm$\vk{P}$} \{#1\} }

\newcommand{\pb}[1]{\mbox{\rm$\vk{P}$}\Bigl \{#1 \Bigr \}}

\newcommand{\R}{\!I\!\!R}

\newcommand{\inr}{\in \R}

\newcommand{\ldot}{,\ldots,}

\newcommand{\limit}[1]{\lim_{#1 \to   \infty}}

\newcommand{\todis}{\stackrel{d}{\to}}

\newcommand{\ntoi}{n \to \infty }
\newcommand{\equaldis}{\stackrel{d}{=}}

\newcommand{\BQN}{\begin{eqnarray}}
\newcommand{\EQN}{\end{eqnarray}}
\newcommand{\BQNY}{\begin{eqnarray*}}
\newcommand{\EQNY}{\end{eqnarray*}}
\newcommand{\BS}{\begin{sat}}
\newcommand{\ES}{\end{sat}}
\newcommand{\BL}{\begin{lem}}
\newcommand{\EL}{\end{lem}}
\newcommand{\BT}{\begin{theo}}
\newcommand{\ET}{\end{theo}}
\newcommand{\BK}{\begin{korr}}
\newcommand{\EK}{\end{korr}}
\newcommand{\BD}{\begin{de}}
\newcommand{\ED}{\end{de}}
\newcommand{\BIT}{\begin{itemize}}
\newcommand{\EIT}{\end{itemize}}
\newcommand{\BDI}{\begin{description}}
\newcommand{\EDI}{\end{description}}

\newcommand{\QED}{\hfill $\Box$}

\newcommand{\IF}{\infty}
\def\kal#1{{\cal{ #1}}}
\def\fracl#1#2{\biggr( \frac{#1}{#2} \biggl) }

\newcommand{\prooftheo}[1]{ \textsc{Proof of Theorem} \ref{#1} }

\newcommand{\BEX}{\begin{exxa}}
\newcommand{\EEX}{\end{exxa}}

\def\uunJ{\vk{u}_{n,J}}

\def\NuunJ{\normpAunJ}

\def\ayk{\normpAunJ}
\def\ayk{\chi_n}

\def\konstRR{h_n}

\def\X{\vk{X}}
\def\Y{\vk{Y}}

\def\YNJU{\vk{Z}_{n,I}}
\def\XIabuVN{( \X_I - \AIJ \AJJM \uunJ   ) \Bigl  \lvert \X_{\barI} = \vk{u}_{n,J}}

\def\barI{J}

\def\omgF{{x_F}}
\def\toomgF{ \uparrow \omgF}

\def\VVI{\vk{U}_I}
\def\VVJ{\vk{U}_{J}}
\def\VVV{\vk{U}}

\def\XI{\vk{X}_I}
\def\XJ{\vk{X}_J}
\def\SSI{\vk{S}_I}
\def\SSJ{\vk{S}_J}
\def\SSS{\vk{S}}
\def\AII{A_{II}}

\def\AJJM{A_{JJ}^{-1}}
\def\AIJ{A_{IJ}}
\def\AJI{A_{JI}}
\def\aaJ{ \vk{a}_J}

\def\njd{\{1 \ldot d\}}

\newcommand{\normp}[1]{\lVert #1 \rVert}

\def\normpA#1{{\norm{#1}}_{A}}
\def\normpAunJ{{\norm{\vk{u}_{n,J}}}_{A}}
\def\Wab{W}

\def\bag{\kal{R}_{\alpha,\gamma} }
\def\pst{ p^* }
\def\normAJ{\normpA{\aaJ }}

\def\rabtau{R_{F,g, \normAJ}}
\def\normUNJ{\norm{\uunJ}_{A}}

\def\njk{\{1 \ldot d\}}

\def\QABF{Q_{ F,g,\tau}}
\def\absXD{X_d}

\def\LPS{$W_p$ scale mixture }
\def\LPSo{$W_p$ scale mixture}

\begin{document}

%

\centerline{\Large
Conditional limits of $W_p$ scale mixture distributions}

       \vskip 1 cm

        \centerline{\large Enkelejd Hashorva}
        \vskip 0.8 cm
        \centerline{\textsl{Department of Mathematical  Statistics and Actuarial Science}}
        \centerline{\textsl{University of Bern, Sidlerstrasse 5}}
        \centerline{\textsl{CH-3012 Bern, Switzerland}}
       \centerline{\textsl{enkelejd.hashorva@stat.unibe.ch }}
\today{}
        \vskip 1.4 cm

 {\bf Abstract:} In this paper we introduce the class of $W_p$ scale mixture random vectors with a particular radial decomposition and a independent splitting property specified by some random variable $W_p, p\in (0,\IF)$. We derive several conditional limit results assuming that the distribution of the random radius is in the max-domain of attraction of a univariate extreme value distribution and $W_p$ has a certain tail asymptotic behaviour. As an application we obtain the joint asymptotic distribution of concomitants of order statics considering  certain bivariate \LPS random samples.



\section{Introduction}
Let $(X,Y)$ be a bivariate spherical random vector with stochastic representation
\BQN\label{eq:de}
(X,Y)\equaldis (R U_1, R U_2),
\EQN
where $R$ is a positive random variable with distribution function $F$ being independent of the bivariate random vector $(U_1,U_2)$ which uniformly distributed on the unit circle of $\R^2$ ($\equaldis$ stands for equality of distribution functions). The radial decomposition \eqref{eq:de}
makes the spherical vectors a very tractable class with well-known properties (see e.g. Kotz (2000)).  An  interesting distributional result is that if we condition on one component and consider, say $Z_u^*:=(Y \lvert X=u)$, then the distribution function of $Z_u^*$ can be easily specified in terms of $F$. A canonical example is when $R^2\equaldis X^2,$ with $X$ and $Y$ two independent standard Gaussian random variables
(with mean 0 and variance 1). In this model $Z_u^*, u \inr$ is a standard Gaussian random variable. Remarkably,
even though the random variable $Z_u:= (Y \lvert X>u)$ is not a Gaussian one,  when $u$ tends to $\IF$ it can be approximated (in distribution) by a standard Gaussian random variable.
In fact, interesting enough, both $Z_u$ and $Z_u^*$ can be approximated in distribution by a Gaussian random variable, provided that the
distribution function $F$  is in the Gumbel max-domain of attraction (see Berman (1983), Hashorva (2006)).
Two other important instances
are when $F$ is in the \FRE and Weibull max-domain of attraction (see Berman (1992)), since again $Z_u,$ and $Z_u^*$ can be approximated by some random variables with known distribution (not dependent on $F$).


The radial decomposition \eqref{eq:de} suggests that conditional limit results could be derived without imposing specific distributional assumptions on the random vector $(U_1,U_2)$. The first attempt in this direction is made in Hashorva (2009) where the approximating distribution is in general no longer Gaussian but a polar Kotz distribution.\\
The main interest related to conditional limit results is due to their  important role in the modeling of rare events,
see e.g., Berman (1992),  Joshi and Nagaraja (1995), Ledford and Twan (1998), Abdous et al.\ (2005), Heffernan and Resnick (2007),
Balkema and Embrechts (2007), or  Resnick (2007).\\
  Recent theoretical developments with interesting statistical applications concerning conditional extreme value models are presented (bivariate setup) in the recent deep contributions Heffernan and Resnick (2007) and Das and Resnick (2008a,b).

In this paper we are interested in multivariate distributions with radial decomposition for which conditional limit results such as the Kotz approximation hold under general asymptotic settings. Borrowing the idea of Hashorva (2008a) where beta-independent random vectors are discussed, we introduce in this paper a new class of random vectors
\COM{
Such a random vector $\vk{X}\inr^d,d\ge 2$ possesses the stochastic representation
\BQN\label{1}
\X \equaldis A R \vk{U},
\EQN
with matrix $A\inr^{d\times d}$ and positive associated random radius $R$ with distribution function $F$ being independent of $\vk{U}:=(U_1 \ldot U_d)^\top$, where $\abs{U_1}^p+ \cdots + \abs{U_d}^p=1, p\in (0,\IF)$ (here $^\top$ denotes the transpose sign).
}
with prominent members the Dirichlet and the elliptical random vectors.

Without going to mathematical details we briefly state the main contributions of this paper:\\
a) First the class of \LPS random vectors is introduced; for this model we show a tractable stochastic representation of conditional random vectors
(see \netheo{prop:1} below). This representation, which is of some interest on its own, is the key to our asymptotic results.\\
b) Under certain asymptotic restrictions on the class of \LPS random vectors we  obtain several conditional limit results
 extending those presented in  Hashorva (2008a) for the class of beta-independent random vectors.\\
%
c) Applying the Kotz approximation of \LPS random vectors we derive the joint asymptotic distribution of the concomitants of order statistics (bivariate setup).
It is interesting that for our model, which includes as the special case the  bivariate elliptical one,
the concomintants of order statistics are asymptotically independent.

Three other applications of our results (not developed here) concern i) joint tail asymptotics of bivariate \LPS random vectors (see Hashorva (2008a)),  ii) the asymptotic independence and asymptotic behaviour of extremes of \LPS random sequences, and iii) the estimation of conditional distributions and conditional quantile functions in multivariate \LPS models (see Abdous et al.\ (2008), Hashorva (2008a,2009)).

Organisation of the paper: The main results are given in Section 3. In Section 4 we present the application to concomitants of order statistics. Proofs of the stated results are relegated to Section 5.

\section{Preliminaries}
We introduce first some standard notation. Throughout  this paper $I$ is a non-empty index subset of $\{1 \ldot d\}, d\ge 2,$
with $m:=\abs{I}< d$ elements. For given $\x=(x_1 \ldot x_d)^\top \inr^d$ we define the subvectors  $\x_I:=(x_i,i\in I)^\top, \x_J:=(x_i,i\in J)^\top$, with $J:= \{1 \ldot d\}\setminus I$ (here $^\top$ denotes the transpose sign). Similarly, given   a matrix $A\inr^{d\times d}$ we define $A_{II},A_{IJ}, A_{JI}, A_{JJ}$  for its  submatrices obtained by keeping the rows and columns with the indices in the corresponding index sets. If $(A_{JJ})^{-1}$ exists we write instead $A_{JJ}^{-1}$.
Next, let $\y=(y_1 \ldot y_d)^ \top $ be another vector in $\R^d$ and let $\norm{\cdot}_J$ be a norm in $\R^{d-m}$. We define
$$ \x+\y:=(x_1+y_1\ldot x_d+y_d)^\top ,  $$
$$ \x>\y, \text{ if } x_i>y_i,\quad  \forall\, i=1 \ldot d, $$
$$\x\ge \y, \text{ if } x_i\ge y_i, \quad \forall\,  i=1 \ldot d,$$
$$\x\not= \y, \text{ if for some $i\le d$  }  x_i\not= y_i, $$
$$\a \x:=(a_1x_1\ldot a_dx_d)^\top,\quad c\x:=(cx_1\ldot cx_d)^\top, \quad  \a\inr^d, c\inr,$$
$$ \norm{\x_J}_A:= \normp{A_{JJ}^{-1} \x_J}_J,\quad A\inr^{d\times d}. $$
In order to simplify the notation, we write  $\x_I^\top$ instead of $(\x_I)^\top$, respectively.\\
Given  $\vk{Z}$   a random vector with distribution function $G$ and density function $g$   we write alternatively $\vk{Z} \sim G$, and $\vk{Z} \backsimeq g$, respectively. When $G$ is a univariate distribution we denote by  $\overline G$ its survivor function and by $x_G$ its upper endpoint.
In the following $Z \sim Beta(a,b)$ or $Z\sim Gamma(a,b)$ mean that $Z$ is Beta or Gamma distributed with positive parameters $a,b$, respectively.
The corresponding density functions are
$$x^{a-1}(1- x)^{b-1}\Gamma(a+b)/(\Gamma(a)\Gamma(b)), x\in (0,1), \quad  \text{ and  }x^{a-1} \exp(- b x) b^a/\Gamma(a), x\in (0,\IF) ,$$ with $\Gamma(\cdot)$ the Gamma function.\\

Let  $I, p\in (0,\IF)$ be given and throughout this
paper $\vk{U}$ is  a random vector in $\R^d,d\ge 2$ with specific properties defined in terms
$\norm{\cdot}_I,\norm{\cdot}_J$ two given norms in $\R^m$ and $\R^{d-m}$, respectively. Explicitly, we suppose that
$$\normp{\VVI}_I=\normp{\VVJ}_J=1$$
 almost surely  and $\VVI$ is independent of $\VVJ$. In our definition below, the properties of $\vk{U}$ are related to $p$ and a random variable $W_p\in [0,1]$.

We arrive at the following definition:

\BD
 Let $R\sim F, W_p \sim G$ be two  random variables such that $R>0, W_p\in [0,1]$ almost surely. We define a \LPS random vector $\vk{X}$ in $\R^d, d\ge 2$ via the stochastic representation
\BQN \label{rep:elli0:b}
 \X \equaldis A\SSS, \quad   \SSI \equaldis R \Wab\VVI , \quad \SSJ \equaldis  R W_p\VVJ, \quad W:=(1- W_p^p)^{1/p},
\EQN
where  $A\in \R^{d\times d}$, and $R, \Wab, \VVI, \VVJ$ are mutually independent.
\ED
Clearly, the distribution function of $\X$ with stochastic representation \eqref{rep:elli0:b}
is determined by $A,F,G,I,p$ and $\VVI, \VVJ$. In the following
we refer to $\X$ as a \LPS random vector with parameters $A,F,G,I,p$, or simply as a \LPS random vector. If $G$ possesses a density function $g$ (which we assume in the following) we mention $g$ instead of $G$. Both $\VVI$ and $\VVJ$ are also important, however since we deal with the same $\vk{U}$ we do not consider
these random vectors as  further parameters in our definition.\\
 In the special case that $W_p^\delta ,\delta\in (0,\IF)$ is a Beta distributed random variable and 
\BQNY
 \norm{\x_I}_I:= \Bigl(\sum_{i\in I} \abs{x_i}^p\Bigr)^{1/p}, \quad  \norm{\x_J}_J:= \Bigl(\sum_{i\in J} \abs{x_i}^p\Bigr)^{1/p}, \quad \x \inr^d,
 \EQNY
$\X$ is referred to as a beta-independent random vector (see Hashorva (2008a)).

Any spherical random vector $\vk{X}$ in $\R^d$ with positive associated random radius satisfies \eqref{rep:elli0:b}  with $\Wab^2\sim Beta(m/2,(d-m)/2)$ and $\vk{U}_I,\vk{U}_J$ two independent random vectors being uniformly distributed on the unit spheres of $\R^{m}$ and $\R^{d-m}$, respectively ($\norm{\cdot}_I$ and $\norm{\cdot}_J$ are the $L_2$-norms in the corresponding spaces). Also Dirichlet and beta-independent random
vectors belong to the class of  \LPS random vectors.

We impose further the following assumption on $\X$ needed in the definition of the conditional random vector $\X_I \lvert \X_J$.

{\it A1.} If $\X$ is a random vector with stochastic representation {\eqref{rep:elli0:b} and $K\subset J$ where $J\setminus K$ has only one element,
 then $\vk{U}_K,$ possesses a positive density function defined for all $\vk{u}_K\inr^{ \abs{K}}$ with $\vk{u}\inr^d$ such that $\normp{\vk{u}_J}= 1$. When $J$ has only one element, then we suppose that $\pk{\vk{U}_J=1}>0$.  Further, we suppose that  $A_{JJ}^{-1}$ exits.\\

In the bivariate setup $d=2$ the random vector $(X,Y)$ defined by
\BQN\label{eq:exa:biv}
 X= R I_1 W, \quad Y= \rho R I_1 W+ (1- \rho^p)^{1/p}  R I_2 W_p, \quad p\in (0,\IF), \quad \rho \in (-1,1),
 \EQN
where $R>0, W:=(1- W^p_p)^{1/p}  \in [0,1]$ almost surely and $I_1,I_2$ assume values in $\{-1,1\}$ with $\pk{I_1=1}\pk{I_2=1}\in (0,1]$ is
a \LPS random vector, provided that $I_1,I_2,R,W$ are mutually independent.


Let  $g:(0,1) \to (0,\IF)$ be the  density function of $W_p $,
and define the distribution function $\QABF$ by
\BQN\label{eq:dfQ}
\QABF(z):=1- \frac{\int_{(\tau^p+z^p)^{1/p}}^\omgF g(\tau/r) \frac{1}{r} \, d F(r)}
{\int_{\tau}^\omgF g(\tau/r) \frac{1}{r}\, d F(r)} ,
\quad \forall z \in (0,(\omgF^p-\tau^p)^{1/p}), \quad \tau \in (0, x_F).
\EQN
We have the following result:
\BT\label{prop:1}
Let $\X$ be a \LPS random vector in $\R^d,d\ge 2$ with representation \eqref{rep:elli0:b} and parameters $A,F,g,I,p$.
Assume that $\AJI$ has all entries equal 0, and Assumption A1 is satisfied.
If  $\a\inr^d$ is such that $F(\normAJ) \in (0,1)$, then we have the stochastic representation
\BQN
\label{eq:res1}
\Bigl( \XI \Bigl \lvert  \XJ = \aaJ \Bigr) &\equaldis &
\AII \rabtau \VVI+ \AIJ \AJJM\aaJ,
\EQN
with $\VVI$ independent of the positive random variable $\rabtau $ which has distribution function $Q_{F,g,\normAJ}$.  \ET

{\bf Remarks:} {\it
(a) If $F$ in \netheo{prop:1} possesses a distribution function $f$, then also $Q_{F,g,\normAJ}$
possesses a density function given by (set $c:=\normAJ \in (0,x_F))$
\BQN\label{eq:density:R}
q_{f,g, c}(z):= \frac{z^{p-1} }{c^p+z^p} \frac{g(c/(c^p+z^p)^{1/p})f((c^p+z^p)^{1/p})}{{\int_{c}^\omgF
g(c/s) \frac{1}{s}f(s)\, d s}  }, \quad \forall z \in (0,(\omgF^p-c^p)^{1/p}).
\EQN
(b) When $\X$ is an elliptical random vector, then \eqref{eq:res1}  holds for any $I\subset \njk$ without the assumption that $\AJI$ has all entries equal 0.\\
(c) Examples of norms appearing in the definition of the random vector $\vk{U}$ are 
\BQN\label{eq:Nq}
 \norm{\x_I}_I:= \Bigl(\sum_{i\in I} \abs{x_i}^{q_1}\Bigr)^{1/q_1}, \quad  \norm{\x_J}_J:= \Bigl(\sum_{i\in J} \abs{x_i}^{q_2}\Bigr)^{1/q_2}, \quad \x \inr^d, \quad q_1,q_2\in [1,\IF).
 \EQN
We note in passing that our results can be stated also when $q_1,q_2\in (0,1)$.
}

The main asymptotic condition imposed on the distribution function $F$ is that it belongs
to the max-domain of attraction of a univariate extreme value distribution function $H$. Explicitly, we suppose that
\BQN
\label{eq:LL}
 \limit{n}\sup_{x\inr} \Abs{F^n(a_nx+ b_n)- H(x)}&=& 0
 \EQN
holds for some $a_n>0,b_n\inr,n\ge 1$. The distribution function $H$ is either the unit Gumbel distribution $\Lambda(x):= \exp(-\exp(-x)),x\inr$,
the unit Weibull distribution $\Psi_\gamma(x):= \exp(- \abs{x}^\gamma), \gamma \in (0,\IF), x\in (- \IF,0)$, or the
unit \FRE distribution $\Phi_\gamma(x):= \exp(- x^{-\gamma}), \gamma\in (0,\IF), x\in (0,\IF)$. See
Reiss (1989), Embrechts et al.\ (1997), Falk et al.\ (2004), De Haan and Ferreira (2006), or
Resnick (2008) for details on univariate extreme value distributions and max-domains of attraction.

\def\ralp{\kal{R}_{\alpha}}
\def\Yy{\kal{E}}

\section{Kotz Approximation}
Consider  $\X$ a \LPS random vector in $\R^d,d\ge 2$  with representation \eqref{rep:elli0:b}
 and parameters $A,F,g,I,p$, and let  $\vk{u}_n,n\ge 1$ be a sequence of constants in $\R^d$  satisfying $\normp{\vk{u}_{n,J}}_A \in (0,\omgF),n\ge 1$. Next, we introduce two sequences of random vectors $\vk{Z}_{n},n\ge 1,$ and $\vk{Z}_{n}^*, n\ge 1$  defined in the same probability space satisfying
\BQN \label{eq:YNJ}
\YNJU&\equaldis &\XIabuVN, \quad \vk{Z}_{n}^*\equaldis  \Bigl( (\X_I -\vk{u}_{n,J}
A_{IJ}), (\X_J- \vk{u}_{n,J})\Bigr)\Bigl \lvert  \X_J> \vk{u}_{n,J}, \quad n\ge 1.
\EQN
For notational simplicity we write $\x_{n,K}$ instead of  $(\x_n)_K$ for some $\x_n\inr^d$ and $K \subset \{1 \ldot d\}.$ \\
Our main concern in this section is the asymptotic approximation  of these sequences when  $\normp{\vk{u}_{n,J}}_A \to  \omgF$ assuming further
that $F$ is in the Gumbel max-domain of attraction satisfying \eqref{eq:LL} with $H=\Lambda$.
The latter assumption (henceforth abbreviated $F\in MDA(\Lambda, w)$ is equivalent with (see e.g., Embrechts et al.\ (1997))
\BQN \label{eq:rdfd}
\lim_{t \uparrow \omgF} \frac{\overline{F}(t+x/w(t))}
{\overline{F}(t)} &=& \exp(-x),\quad \forall x\inr,
\EQN
where $w$ is some positive scaling function. If $\X$ is an elliptical random vector i.e., $\vk{U}$ is
uniformly distributed on the unit sphere of $\R^d$ (with respect to $L_2$-norm), then in view of Hashorva (2006) both conditional random vectors
$\YNJU$ and  $\vk{Z}_{n,I}^*$ can be approximated in distribution by a Gaussian random vector, provided that $\normp{\vk{u}_{n, J}}_A \to x_F$.\\
We note in passing that the Gaussian approximation for bivariate elliptical random vectors is first obtained in full generality in Berman (1983).
In Hashorva (2008a) it is shown that the limiting random vector
 is a Kotz Type I polar random vector if $\X$ is a beta-independent random vector. In our definition $\Y \in \R^k,k\ge 2$ is referred to as a Kotz Type I scale mixture random vector if $\vk{Y}\equaldis B R \vk{V}$, where $B\inr^{k\times k}$ and
$R^q \sim Gamma(\alpha, \beta), \alpha, \beta, q \in (0,\IF)$. Further $R$ is positive and independent
of $\vk{V}\in\R^k$ and for some norm $\norm{\cdot}_{k}$ in $\R^k$ we have  $\normp{\vk{V}}_k=1$ almost surely. When $\normp{\cdot}_k$ is 
the Euclidian norm, $\X$ is referred to as a Kotz Type I polar random vector. 

In this section we show that Kotz approximation of conditional random vectors $  \YNJU$ and $\vk{Z}_{n,I}^* $ holds also for the
general settings of \LPS random vectors. Instead of some distributional assumptions on $\vk{U}$, we suppose only that
the density function $g$ of $W_p$ satisfies
\BQN
\label{eq:g:alp}
\limit{u} \frac{ g(1- x/u)}{g(1- 1/u)}= x^{\alpha-1}, \quad \forall x>0
\EQN
for some $\alpha\in (0,\IF)$, i.e., $g$ is regularly varying at 1 with index $\alpha-1$. As will be shown below the parameter $\alpha$ together with $p$ determines the conditional limit distribution.

We state next the main result of this section.
\BT \label{eq:mainT}
Let  $I$ be an index set of $\njd, d\ge 2$, and let $\X$ be a \LPS  random vector in $\R^d$
with representation \eqref{rep:elli0:b} and parameters $A,F,g,I,p$ satisfying Assumption A1.
Suppose that \eqref{eq:g:alp} holds with $\alpha\in (0,\IF)$ and for any $\ve\in (0,1)$
\BQN \label{eq_cond:GUM}
g(x)\le c_{\ve} x^{ \gamma_\ve}, \quad \forall x \in (0,\ve)
\EQN
is valid with $c_\ve\in (0,\IF), \gamma_\ve\inr$. \resub{Assume further that $F\in MDA(\Lambda, w)$ and 
$\AJI$ has all entries equal 0. }

a) Let $\vk{u}_n,n\ge 1$ be constants in $\R^d$ such that $F(\normp{\vk{u}_{n, J}}_A)\in (0,1),n\ge 1$
and $\limit{n}\normp{\vk{u}_{n, J}}_A=\omgF$. Then we have the convergence in distribution
\BQN\label{eq:main:1}
\konstRR \YNJU &\todis & \AII \ralp \VVI, \quad n\to \infty,
\EQN
where $\konstRR :=(\normp{\vk{u}_{n, J}}_A^{1-p}w(\normp{\vk{u}_{n, J}}_A) )^{1/p}, n>1$, and $\ralp$ is a positive random variable independent of $\VVV_I$ such that $\ralp^{p}\sim Gamma(\alpha,  1/p)$.\\
b) Let $\vk{u}_n:=u_n(1\ldot 1)^\top \inr^d, u_n\in (0,\omgF), n\ge 1$ be such that $\limit{n} u_n= \omgF$.
If  $I= \{1 \ldot d-1\}, J=\{d\},A_{JJ}=1$, and $\AIJ$ has all entries equal 0 if $p\in (0,1)$, then we have
 (set $1_p:=1$ if $p=1,$ and $1_p:=0$ otherwise)
\BQN\label{eq:main2:2P}
\Bigl( h_n \vk{Z}_{n,I}^* , w(u_n)\vk{Z}_{n,J}^* \Bigr) &\todis &
\Bigl(\AII \ralp  \VVI+ \Yy 1_p\AIJ, \Yy\Bigr), \quad \ntoi,
\EQN
with $\Yy\sim Gamma(1,1)$ independent of $(\ralp, \VVI)$.
\ET

\begin{remark}
(a) \cre{
The random vector $\vk{Y}_I:= \AII \ralp \VVI$ appearing in \eqref{eq:main:1} is a Kotz Type I scale mixture random vector.
Therefore we refer to the distribution approximation in \eqref{eq:main:1} as the Kotz approximation.} It
is well-known (see e.g., Kotz et al.\ (2000)) that $\vk{Y}_I$ is a Gaussian random vector in
$\R^{m}, m:=\abs{I}$ with covariance matrix $A_{II} A_{II}^\top$, provided that $p=2,\alpha=m/2,$
and $\VVI$ is uniformly distributed on the unit sphere of $\R^{m}$ (with respect to $L_2$-norm).\\
(b) The Kotz approximation above is stated in terms of convergence of distribution functions. It is of some interest to strengthen this  to convergence of the corresponding density functions, which we refer to as the strong Kotz approximation.
The random vector $\vk{Z}_{n,I}, \vk{Z}_{n,I}^*, n\ge 1 $ possess a density function (recall \eqref{eq:density:R}) if both $R$ and $\VVI$ possess a positive density function. If $R\backsimeq f$ such that
\BQN\label{eq:local:gumbel}
\lim_{u\uparrow x_F} \frac{f(u+x/w(u))}{f(u)} &=&\exp(-x), \quad \forall x\inr
\EQN
holds with some positive scaling function $w$, then for $\X$ a $L_p$ Dirichlet random vector (see Hashorva and Kotz (2009)) the convergence in \eqref{eq:main2:2P} can be strengthened to the strong Kotz approximation. With similar arguments as in the aforementioned paper utilising further \eqref{eq:density:R}, it follows that both \eqref{eq:main:1} and \eqref{eq:main2:2P} can be strengthened to the local uniform convergence of the corresponding density functions, provided that \eqref{eq:local:gumbel} is satisfied and $\VVI$ possesses a positive density function.
\end{remark}

We present next two illustrating examples.

{\bf Example 1.} [Kotz Type III \LPSo] We refer to  $\X$ in $\R^d,d\ge 2$ as a Kotz Type III \LPS random vector if it has stochastic representation \eqref{rep:elli0:b} (for some given index set $I$)
with $A  \in  \R^{d\times d}$ such that $A_{JJ}$ is non-singular, and
\BQN\label{eq:kotz:Fk}
\overline{F}(u) &=&  (1+o(1))K  u^{N}\exp(-r u^\delta), \quad K>0,\delta>0, N\inr, \quad u\to \infty.
 \EQN
Assume further that $W_p$ possesses the positive density function $g$ which is bounded in $[0,1]$ and for any $\ve \in (0,1)$
$$ g(u)=c(u)(1- u)^{\alpha-1}, \quad \forall u\in (\ve, 1), \quad \alpha\in (0,\IF),$$
with $c(u)$ some positive measurable function such that $\lim_{u\uparrow 1} c(u)=c\in (0,\IF)$.
Clearly $g$ satisfies the assumptions of \netheo{eq:mainT}. Since $F$ is in the Gumbel max-domain of attraction with scaling function $w(s)=r \delta s^{\delta -1}, s>0$ the aforementioned theorem implies
for any sequence $\vk{u}_n,n\ge 1$ satisfying $\limit{n} \normp{\vk{u}_{n, J}}_A= \IF$
\BQNY
(r \delta )^{1/p}\normp{\vk{u}_{n, J}}_A^{\delta/p - 1} \YNJU &\todis & \AII \ralp \VVI, \quad n\to \infty,
\EQNY
where $\ralp$ is a positive random variable independent of $\VVV_I$ such that $\ralp^{p}\sim Gamma(\alpha,  1/p), \ralp\in (0,\IF)$.

{\bf Example 2.} [$F$ with finite upper endpoint] Let $\X= A R \vk{U}$ be a \LPS random vector in $\R^d,d\ge 2$. Suppose that $R\sim F$ with
$$ \overline{F}(u)= (1+o(1))c_1 \exp(- c_2(x_F- u)^{-\lambda}) , \quad c_1,c_2,\lambda \in (0,\IF) \quad u \uparrow x_F\in (0,\IF).$$
It follows easily that $F\in MDA(\Lambda, w)$ with $w(s)= c_2 \lambda (x_F- s)^{- \lambda-1}, s\in (0,1)$.
Consequently, if $A, I,\vk{U}, \vk{u}_n,n\ge 1$ are such that the assumptions of \netheo{eq:mainT} are satisfied, then we have the approximation
\BQNY
\frac{(c_2 \lambda)^{1/p} x_F^{1/p- 1}}{ (x_F- \normp{\vk{u}_{n,J}})^{(\lambda+1)/p}} \YNJU &\todis & \AII \ralp \VVI, \quad n\to \infty,
\EQNY
where $\limit{n} \normp{\vk{u}_{n,J}}= x_F$ and $\ralp, \VVV_I$ are as in Example 1.

\def\coEHJ#1{#1}
\def\barAL{\overline{\alpha}}
\def\barALI{\overline{\alpha}_I}
\def\barALJ{\overline{\alpha}_J}
\def\barALK{\barAL- \barALI-\barALJ}

\def\konstRRb{h_n}
\section{Regularly Varying $\overline{F}$}
In this section we deal with distribution functions $F$ in the Weibull or \FRE max-domain of attraction. Specifically, in the former case
\BQN\label{eq:mean:psi}
 \limit{u} \frac{\overline{F}(1- x/u)}{\overline{F}1-1/u)}& =& x^{\gamma}, \quad \forall x>0, \quad \gamma \in (0,\IF)
\EQN
is valid for some $\gamma \in (0,\IF)$. Condition \eqref{eq:mean:psi} is equivalent with  \eqref{eq:LL} where $H=\Psi_\gamma$ is the unit Weibull distribution. When $F$ is in the max-domain of attraction of the \FRE distribution $\Phi_\gamma$, then we have
\BQN\label{eq:mean:FRE}
 \limit{u} \frac{\overline{F}(xu)}{\overline{F}(u)}&=& x^{-\gamma}, \quad \forall x>0.
 \EQN
Under these assumptions, where the survivor function $\overline{F}$ is regularly varying at the upper endpoint, the approximation of the conditional distribution of \LPS random vectors can be carried out as in the case of beta-independent random vectors.\\
It is interesting that when \eqref{eq:g:alp} and \eqref{eq:mean:psi} hold, then the conditional
limit distributions are completely specified by $\alpha, p$ and $\gamma$.
For the \FRE case we do not to impose any asymptotic assumptions on $g$, therefore the conditional limit distribution depends only on the regularly varying index $\gamma$. Both $\alpha$ and $p$ do not appear in the asymptotics. 

We state next the main result of this section.

\def\normAAJ{\normpA{\vk{a}_{J}}}
\def\RRabp{\kal{R}_{\normAAJ,\alpha,\beta,\gamma,p,\delta}}

\BT \label{theo:weib:1}
Let $\X, \vk{u}_n,\YNJU,n\ge 1,$ be as in \netheo{eq:mainT}, and let $\bag, -\Yy$ be two independent positive random variables such that $ \bag^{p}\sim Beta(\alpha, \gamma)$, and $\Yy\sim 1 - \abs{x}^{\alpha+\gamma}, x\in (-1,0)$  being both independent of $\VVI$.
Assume that $F$ satisfies \eqref{eq:mean:psi}, and  let 
$a_n,n\ge 1$ be positive constants such that $\limit{n}a_n=0$.\\
a) If for all large $n$ we have $ \NuunJ=1- a_n$, then (set $\konstRRb:=(pa_n)^{-1/p}, n\ge 1$)
\BQN\label{eq:main:1:2}
\konstRRb \YNJU &\todis & \AII \bag  \VVI, \quad n\to \infty.
\EQN
b) If $I=\{1 \ldot d-1\}, A_{JJ}=1$, and $\AIJ$ has all entries equal 0 if $p\in (0,1)$, then for any  $\x\in \R^ d$
with
$x_d\in (-\IF, 0)$
\BQN\label{eq:main:1:3}
\lefteqn{\limit{n}\pb{ \konstRRb (\X_I-  \AIJ  ) \le \x_I, (X_d-1)/a_n \le x_d \bigl \lvert X_d > 1- a_n  }}\notag\\
&=& \pk{\abs{\Yy}^{1/p}(\AII \bag \VVI - 1_p \AIJ  )\le \x_I, \Yy \le x_d}
\EQN
is valid. \\
c) If \eqref{eq:mean:FRE} holds, 
then for any $\a\inr^d$ such that $\normAAJ >0$ we have (set $\vk{u}_n:= \vk{a}/a_n, n\ge 1 $)
\BQN\label{eq:main:1:5}
 a_n  \vk{Z}_{n,I} &\todis & \AII \kal{R} \VVI, \quad n\to \infty,
\EQN
with $\kal{R}>0$ independent of $\VVI$ satisfying $ \kal{R}\sim Q_{M, g, \normAJ},$
where $M(s):= 1- \normAAJ^{\gamma} s^{-\gamma}, \forall s\ge \normAAJ.$
\ET

\begin{remark}
(a) The convergence in \eqref{eq:main:1:2} and \eqref{eq:main:1:5} is stated in Hashorva et al.\ (2007)
for $L_p$ Dirichlet random vectors, and Hashorva (2008a) for beta-independent random vectors. \\
(b) Under von Misses conditions on the density function $f$ of $F$ the above asymptotics can be strengthened to local
uniform convergence of the corresponding density functions.\\
(c) Condition \eqref{eq:mean:psi} is satisfied for instance for the Beta distribution, whereas condition
\eqref{eq:mean:FRE} is satisfied by distributions $F$ with tail behaviour $\overline{F}(x)= (1+o(1))\lambda  x^{-\gamma}, \lambda ,\gamma \in (0,\IF)$ as $x\to \IF$.\\
(d) When $\X$ is a \LPS random vector, then $\X/c, c\in (0,\IF)$ is also a \LPS random vector. Further  if  $R\sim F$ with $F$ in the max-domain of attraction of $\Psi_\gamma$, then  also $R/c$ has distribution function in the max-domain of attraction of $\Psi_\gamma$. Hence  the extension of our asymptotic results in \eqref{eq:main:1:2} and \eqref{eq:main:1:3} for  $F$ with upper endpoint $x_F\in (0,\IF)\setminus  \{1 \}$ follows easily.
\end{remark}

\section{Asymptotics of Concomitants of Order Statistics}
We present next an application of the Kotz approximation concerning the asymptotic distribution of concomitants of order statistics from a bivariate sample with underlying \LPS  distribution. Let therefore $(X_{j},X_{j}), j=1 \ldot n$ be independent bivariate random vectors with stochastic representation \eqref{eq:exa:biv}. The $i$th concomitant of order statistics $Y_{[i:n]}$ is defined as follows:
if we order the pairs based on the order statistics $X_{1:n} \le \cdots \le X_{n:n}$, then
$Y_{[i:n]}$ is the second component of the pair with first component the $i$th order statistics $X_{i:n}$.\\
The main applications of concomitants are in selection procedures, ranked-set sampling, prediction analysis, and
 inference problems,  see e.g., David and Nagaraja (2003), or Wang (2008) for details.

In an asymptotic context, Nagaraja and David (1994) derive interesting results for the maximum of the concomitants of order statistics
$Y_{n,k}:=\max_{1 \le i\le k}Y_{[n-i+1:n]}, n>1,$ with $k$ a fix integer. Related asymptotic results can be found in Eddy and Gales  (1981), Galambos (1987), Joshi and Nagaraja (1995), and Ledford and Tawn (1998).

In fact,  by Result 1 of Nagaraja and David (1994) and the Kotz approximation the asymptotic distribution of
$Y_{n,k}, n \to \IF$ follows easily. Below we shall investigate the asymptotic behaviour of $(Y_{[n:n]}, \ldot Y_{[n-k+1:n]}), 1 \le k <n$.
Since the joint distribution function of $(Y_{[n:n]}\ldot Y_{ [n-k+1:n]}), 1 \le k< n$ conditioning on the order statistics
$(X_{n:n}=x_1 \ldot X_{n-k+1:n}=x_{k})$  is
\BQN \label{eq:concos}
\prod_{i=1}^k \pk{Y_1 \le y_i \lvert X_1=x_i},  \quad x_i,y_i\in \R
\EQN
we can apply the Kotz approximation developed in Section 3. More specifically, we shall show that $(Y_{[n:n]}, \ldot Y_{[n-k+1:n]})$ are asymptotically independent, provided that the \LPS random vector $(X_1,Y_1)$ satisfies the assumptions of \netheo{eq:mainT}. It is interesting that in this model
the limit distribution of the  concomitants of order statistics depends only on $\alpha$  and $p$.\\

\BT \label{theo:coco} Let $(X_j,Y_j), j=1 \ldot n$ be a random sample of \LPS  random vectors as defined in \eqref{eq:exa:biv}. Assume that $p\in (1,\IF), R \sim F$, with $F$ in the Gumbel max-domain of attraction with some scaling function $w$, and  $W_p\backsimeq g$. Define $A_n,B_n$ by
$$ A_n:= (1- \rho^p)^{1/p} \frac{b_n}{(b_n w(b_n))^{1/p}}, \quad B_n:= \rho b_n, \quad b_n:=H^{-1}(1- 1/n), \quad n>1,
$$
with $H^{-1}$ the inverse of the distribution function of $X_1$. If \eqref{eq:g:alp} holds with $\alpha\in (0,\IF)$, and the positive density $g$ satisfies \eqref{eq_cond:GUM},  then for any $k\ge 1$ we have
\BQN\label{ANBN}
\Biggl( \frac{Y_{[n:n]}- B_n}{A_n} \ldot  \frac{Y_{[n-k+1:n]}- B_n}{A_n}\Biggr)
&\todis & (\eta_1 \ldot \eta_k), \quad n\to \IF,
\EQN
where  $\eta_1 \ldot \eta_k$ are independent random variables being symmetric about 0 satisfying  $\abs{\eta_i}^p \sim \Gamma(\alpha, 1/p), i\le n$.
\ET
In \netheo{theo:coco} we discuss only the case $p\in (1, \IF)$. When $p\in (0,1]$ and $\rho=0$ the same asymptotic result holds.

If $(X_1,Y_1)$ is a standard bivariate Gaussian random vector, then \eqref{ANBN} implies
\BQN\label{eq.notnow}
\Biggl( \frac{Y_{[n:n]}- \rho \sqrt{2 \ln n} }{\sqrt{1- \rho^2} } \ldot  \frac{Y_{[n-k+1:n]}- \rho \sqrt{2 \ln n }}{
\sqrt{1- \rho^2}}\Biggr) &\todis & (\eta_1 \ldot \eta_k), \quad n\to \IF,
\EQN
with $\eta_i, i\le n$ independent standard Gaussian random variables. Hence  for any $k\ge1 $ we obtain (see Nagaraja and David (1994),
Joshi and Nagaraja (1995))
\BQN\label{eq.notnow2}
\frac{Y_{n,k}- \rho \sqrt{2 \ln n} }{\sqrt{1- \rho^2} } &\todis & \max_{1 \le i \le k} \eta_i, \quad n\to \IF.
\EQN
Further, remark that under the assumptions of \netheo{theo:coco}
\BQN \label{eq:strongAA}
\Biggl( \frac{Y_{[n-i+1:n]}- B_n}{A_n}, \frac{Y_{n-i+1:n}- b_n}{a_n}\Biggr) &\todis &(\eta_i, \xi_i), \quad n\to \IF,
\EQN
where $\eta_i $ is independent of $\xi_i, i\ge 1$ and $\xi_i \sim \Lambda'  (-  \ln \Lambda)^i/i!$ with $\Lambda'$ the density of $\Lambda$.\\
\COM{
When $R\backsimeq  f$ the bivariate random vector $(Y_{[n-i+1:n]}, Y_{n-i+1:n})$ possesses
the joint density function
$$ h(y \lvert x) q_{i:n}(x),$$
with $h(y \lvert x)$ the conditional density function of $Y \lvert X=x$ and $q_{i:n}$ the density function of the order statistics $X_{n-i+1:n}$. When $f$ satisfies \eqref{eq:local:gumbel}, then the strong Kotz approximation implies that \eqref{eq:strongAA}
can be strengthened to local uniform convergence of the corresponding density functions.
}

\section{Proofs}
We present next a lemma and then proceed with the proofs.
\BL \label{lem:1}
Let $g:[0,1]\to (0,\IF)$ be a positive measurable function, and let $F$ be a distribution function on $[0,\IF)$ with upper endpoint $x_F\in (0,\IF]$.\\
a) Let $y \in (0,\IF),z\in [y, \IF)$ be given constants. If $F$ is in the max-domain of attraction of $\Phi_\gamma, \gamma\in (0,\IF)$ and
\resub{$g(y/r)\le cr^{\delta},\forall r\ge z, $ with $\delta < \gamma+1$} and $c\in (0,\IF)$, then  we have
\BQN\label{eq:abo}
\int_{uz}^\IF  g( uy/r)\frac{1}{r} \, d F(r) &=&
(1+o(1)) \overline F(u) \gamma \int_{z }^{\IF} g(y/r)r^{-\gamma-2}\, dr , \quad u \to \IF.
\EQN
In the special case $g$ is a density function of a positive random variable $Z$, then \eqref{eq:abo} holds,  provided that
 $\E{Z^{\tau}}\in (0,\IF)$ for some $\tau\in (\gamma, \IF)$.\\
b) Suppose that $F$ satisfies \eqref{eq:mean:psi}. If
\eqref{eq:g:alp} holds with $\alpha\in (0,\IF)$, then for any $\beta\inr $ and $\IF >z> y\ge 0$ we have
\BQN
\int_{1 -u(z-y)}^1 g( (1- uz)/x) x^\beta \, d F(x) &=&
(1+o(1)) \gamma \overline F(1- u)g(1- u) \int_{0 }^{z-y}(z- t)^{\alpha-1} t^{\gamma-1} \, dt , \quad u \downarrow 0.
\EQN
 c) If \resub{$F\in MDA(\Lambda, w)$}, then for any given constant $\beta \in\R$ 
\BQN\label{eq:lem1:Kopie}
\limit{u} \frac{ (u w(u))^\beta \overline F(\mu u)}{\overline F(u)}&=&0, \quad \forall \mu  \in (1,\IF)
\EQN
and moreover, if \eqref{eq:g:alp} holds with $\alpha\in (0,\IF)$ and \eqref{eq_cond:GUM} is satisfied, then for any $z\in [0,\IF)$ we have
\BQN\label{eq:lem1:1}
\frac{1}{\overline F(u)}\int_{ u+z/w(u) }^\omgF g(u/x) x^\beta \, dF(x)&=&(1+o(1))g(1- \frac{1}{u w(u)}) u^{\beta} \int_{z}^ {\infty} s^{\alpha-1} \exp(-s)\, ds, \quad u \toomgF.
\EQN
\EL
\begin{proof} a) Assume for simplicity that $F$ possesses a density function $f$ which is also regularly varying.
The  general case follows applying  Lemma 2 in Kaj et al.\ (2007). By the assumption on $g$ the integral
$I_{y,z}:= \int_{z }^{\IF} g(y/r)\frac{1}{r} r^{-\gamma-1}\, d r$ is finite. Further, the regular variation of $f$ implies
$$ \limit{u} \frac{ f(u)}{ u\overline{F}(u)}= \gamma,$$
hence  applying Karamata's Theorem (see e.g., Resnick (2008))  we may write
\BQNY
\int_{uz}^\IF  g( uy/r)\frac{1}{r} \, d F(r) &=&
\frac{f(u)}{u} \int_{z }^{\IF} g(y/r)\frac{1}{r} \frac{f(ur)}{f(u)} \, d r \\
&=& (1+o(1)) \gamma I_{y,z} \overline F(u) , \quad u \to \IF.
\EQNY
If $g$ is a density function of some positive random variable $Z\in (0,1)$, then $I_{y,z}$ is finite if $\E{Z^{\tau}}< \IF, \tau \in (\gamma, \IF)$, and
the statement can be established with the same arguments as above.\\
b) Transforming the variables for any $u\in (0,1)$ we have
\BQNY
\frac{1}{\overline F(1- u)}\int_{1 +uy-uz}^1 g( (1- uz)/x) x^\beta \, d F(x) &=&
\int_0^{y-z} g( (1- uz)/(1- us)) [1- us]^\beta \, d F(1- us )/\overline F(1- u).
\EQNY
Consequently, by the max-domain of attraction assumption on $F$ and the regular variation of $g$ at 1
\BQNY
\frac{1}{\overline F(1- u)}\int_{0}^{y-z} g( (1- uz)/x) x^\beta \, d F(x) &=&(1+o(1))g(1-u) \gamma
\int_0^{y-z} (z-s)^{\alpha-1} s^{\gamma -1} \, d s, \quad u \downarrow 0,
\EQNY
hence the claim follows.
%

c) Set $F_u(x):= F(u+ x/w(u))/\overline F(u),v(u):=u w(u), u>0,x\inr$. The Gumbel max-domain of attraction assumption on $F$ implies (see e.g.\ Resnick (2008))
\BQN \label{eq:res:IMP}
 \lim_{u\uparrow x_F} v(u)=\IF, \quad \lim_{u\uparrow x_F} w(u)(\omgF- u)=\IF \text{  if   }x_F < \IF
 \EQN
and
$$ \lim_{u \uparrow x_F} F_u(s)- F_u(t)= \exp(-t)- \exp(-s), \quad s,t\inr.$$
Transforming the variables for \resub{$u$ large} 
\BQNY
\frac{1}{\overline F(u)}\int_{ u+z/w(u) }^{x_F} g(u/x) x^\beta \, dF(x)
&=& u^\beta  \int_{ z}^{w(u)(x_F- u)} g(1/[1+x/v(u)]) [1+x/v(u)]^\beta \, dF_u(x)\\
&=& (1+o(1))u^\beta g(1- 1/v(u)) \int_{ z}^{\IF} (1+o(1))x^{\alpha-1} dF_u(x).
\EQNY
\COM{Next, Fatou Lemma yields
$$
\liminf_{u\uparrow x_F} \int_{ z}^{w(u)(x_F- u)} (1+o(1))x^{\alpha-1} dF_u(x)\ge
\int_{ z}^{\IF}x^{\alpha-1} \exp(-x) \, dx.
$$
}
\resub{We consider only the case $x_F=\IF$ and omit the proof when $x_F\in (0,\IF)$ since it can be established with the same arguments.}  For any $\ve>0$ we may write
$$ \int_{ u+z/w(u) }^{\IF} g(u/x) x^\beta \, dF(x)=
 \int_{ u+z/w(u) }^{(1+\ve)u} g(u/x) x^\beta \, dF(x)+  \int_{(1+\ve)u}^\IF g(u/x) x^\beta \, dF(x)=:I_{\ve}(u)+ J_{\ve}(u).$$
As in the proof of Lemma 3.5 in Hashorva (2006) utilising further Potter's upper bound (see De Haan and Ferreira (2006) or Resnick (2008)) for the regularly varying function $g$, for any $\ve>0$ sufficiently small we obtain (recall $\limit{u} v(u)=\IF)$
\BQNY \limit{u}\frac{I_{\ve} (u)}{u^\beta g(1- 1/v(u))} &= &
\limit{u}\int_{ z}^{(1+\ve)v(u)} (1+o(1))x^{\alpha-1} dF_u(x)=
\int_{ z}^\IF x^{\alpha-1} \exp(-x) \, dx.
\EQNY
By Lemma 4.5 in the aforementioned paper
\BQNY
I_{\alpha,u}:=\int_u^\IF(x^2- u^2)^\alpha\, d F(x)=(1+o(1)) \Gamma(\alpha+1) \fracl{2 }{v(u)}^\alpha u^{2 \alpha}\overline{F}(u) , \quad u\to \IF.
\EQNY
Since for any $\xi>1$ we have $\limit{u}\overline{F}(\xi u)/F(u)=0,$ and furthermore
\BQNY
I_{\alpha,u}&\ge & \int_{\xi u}^{2\xi u}(x^2- u^2)^\alpha\, d F(x)\ge (1+o(1)) \xi^2 u ^{2 \alpha} \overline{F}(\xi u) , \quad u\to \IF
\EQNY
\eqref{eq:lem1:Kopie} follows easily. Next, by \eqref{eq_cond:GUM}  and \eqref{eq:lem1:Kopie}
\BQNY
J_{\ve}(u)&\le &  c_\ve u^{\gamma_\ve} \int_{(1+\ve)u}^\IF x^{\beta-\gamma_\ve} \, dF(x)= o(I_{\ve}(u)), \quad u\to   \IF,
\EQNY
thus the result follows. \end{proof}

\def\unN{\tau_n}
\def\wun{w(\unN}
\def\normUNJ{\unN}
\def\ayk{\unN}

\prooftheo{eq:mainT} a) Let $R_n,n\ge 1$ be random variable with survivor function
\BQN\label{eq:ER}
\pk{R_n > s}:=\frac{\int_{(\normUNJ^p+s^p)^{1/p}}^\omgF g(\normUNJ/x)\frac{1}{x} \, d F(x)}
{\int_{\normUNJ}^\omgF g(\normUNJ /x) \frac{1}{x} \, d F(x)} , \quad \forall s \in (0,(\omgF^p-\normUNJ^p)^{1/p}),
\EQN
where  $\unN:=\normp{\vk{u}_{n, J}}_A, n\ge 1.$ By \eqref{eq:res:IMP} and \nelem{lem:1} for any $z>0$ we have
\BQNY
\limit{n}\pk{ \konstRR  R_n >z } &=&
 \limit{n}\frac{  \int_{\normUNJ+ (1+o(1))z^p/(p\wun)) }^\omgF g(\normUNJ/x) \frac{1}{x} \, d F(x)}
 {\int_{\normUNJ}^\omgF g(\normUNJ/x) \frac{1}{x}  \, d F(x)}\\
&=&\frac{\int_{z^p/p}^\infty                                                                                                                                                                                                                                                                                                                                                                                                                                                                                                                                                                                                                                                                                                                                                                                                                                                                                                                                                                                                                                                                                                                                                                                                                                                                                                                                                                                                                                                                                                                                                                                                                                                                                                                                                                                                                                                                                                                                                                                                                                                                                                                                                                                                                                                                                                                                                                                                                                                                                                                                                                                                                                                                                                                                                                                                                                                                                                                                                                                                                                                                                                                                                                                                                                                                                                                                                                                                                                                                                                                          t^{\alpha-1}\exp(-t) \,dt} {\int_{0}^\infty                                                                                                                                                                                                                                                                                                                                                                                                                                                                                                                                                                                                                                                                                                                                                                                                                                                                                                                                                                                                                                                                                                                                                                                                                                                                                                                                                                                                                                                                                                                                                                                                                                                                                                                                                                                                                                                                                                                                                                                                                                                                                                                                                                                                                                                                                                                                                                                                                                                                                                                                                                                                                                                                                                                                                                                                                                                                                                                                                                                                                                                                                                                                                                                                                                                                                                                                                                                                                                                                                                                           t^{\alpha-1}\exp(-t) \,dt}\\
&=&\pk{\ralp> z},
\EQNY
with $\ralp$ a positive random variable satisfying $\ralp^p \sim \Gamma(\alpha, 1/p)$. Since
$\ralp$ is independent of $\VVI$ the first claim follows using further \eqref{eq:res1}.\\
b) As in Hashorva (2007c) it follows that $X_d$ has distribution function $H$ in the Gumbel max-domain of attraction with the scaling function $w$. Our proof below is quite similar to the proof of  Theorem 3.1 in Hashorva (2008a), therefore we omit some details.
Next, define $h_n:=(u_n w(u_n))^{(1- p)/p}w(u_n), n\ge 1,$ and let $\x\inr^d, t\inr$ be given. By \eqref{eq:res:IMP} $ \limit{n} h_n/w(u_n)=0$ holds if  $p>1$, and $ \limit{n} h_n/w(u_n)=1$ when  $p=1$. Since further (see e.g., Resnick (2008))
\BQN\label{eq:uv:b}
 \lim_{u \toomgF} \frac{ w(u+z/w(u))}{w(u)}&=&1
\EQN
uniformly for $z$ in compact sets of $\R$  for any  $p\in (0,\IF) $ and $\AIJ$ with elements equal 0 when $p\in (0,1)$ we have
\BQNY
\limit{n} \pb{h_n(\X_I-  u_n \AIJ  ) \le  \x_I \bigl \lvert \absXD = u_n+ t/w(u_n) }
&=&\pk{\AII \ralp\VVI+ 1_p t\AIJ \le \x_I}
\EQNY
locally uniformly for $t\inr$ with $1_p:=1$ if $p=1$ and $1_p:=0$ otherwise.
Hence along the lines of the proof of Theorem 3.3 in Hashorva (2006) for any  $x_d>0$ we obtain
\BQNY
\lefteqn{\pb{ h_n(\X_I-  u_n \AIJ  ) \le \x_I, w(u_n)(\absXD- u_n) \le x_d \bigl \lvert \absXD> u_n }}\\
 &=&
\int_{0}^{x_d} \pb{h_n(\X_I-  u_n \AIJ  ) \le  \x_I \bigl \lvert
\absXD = u_n+ t/w(u_n) }\, d H(u_n+ t/w(u_n))/\overline H(u_n)  \\
 &\to & \int_0^{x_d}\pk{\AII \ralp\VVI+  1_p t \AIJ \le \x_I} \exp(-t)\ dt , \quad \ntoi,
\EQNY
thus the proof is complete. \QED

\prooftheo{theo:weib:1}  a) Let $ R_n, n\ge 1 $ be positive random variables as in \eqref{eq:ER} corresponding to
$\vk{u}_{n, J}$. Since $\limit{n} a_n=0$, by the assumptions on $F$ and $g$
for any $t\in (0,1)$ we obtain (set $\overline{a}_n:=1- a_n, h_n:=(p a_n)^{-1/p}, n\ge 1$)
\BQNY
\limit{n} \pk{h_n R_n >t } &=&  \limit{n}\frac{  \int_{
[\overline{a}_n^p+p a_n t^p]^{1/p}}^1 g(\overline{a}_n/r )\frac{1}{r} \, d F(r)}
 {\int_{\overline{a}_n}^1g(\overline{a}_n/r )\frac{1}{r} \, d F(r)}\\
&=& \limit{n}\frac{  \int_{ 1 -  a_n(1-  t^p)(1+o(1))}^1 g( \overline{a}_n /r )\frac{1}{r} \, d F(r)}
 {\int_{\overline{a}_n }^1g( \overline{a}_n /r )\frac{1}{r} \, d F(r)}\\
&=&  \frac{  \int_{0}^{1-  t^p} (1-x)^{\alpha-1} x^{\gamma-1}\, d x}
{  \int_{0}^{1} (1-x)^{\alpha-1} x^{\gamma-1}\, d x}= 1- B(t^p,\alpha,\gamma),
\EQNY
where $B(x,\alpha,\gamma),x\in (0,1)$ is the Beta distribution function with positive parameters
$\alpha,\gamma$ implying thus
\BQNY
h_n R_n &\todis & \bag,  \quad n\to \infty,
\EQNY
with $\bag >0$ such that  $\bag^p  \sim Beta(\alpha,\gamma)$ being independent of $\VVI$. Consequently,
\BQNY
h_n  \AII R_n  \VVI &\todis & \AII  \bag\VVI, \quad n\to \infty, \EQNY
 hence the first claim follows.

 b) In view of Theorem 3.1 in Hashorva (2008c) the random variable $X_d$ has distribution function $H$
 in the max-domain of attraction of $\Psi_{\alpha +\gamma}$.
 If $\vk{u}_{n,J}= 1-a_nt, n\ge 1, t>0$, then for any $\x \in \R^d$ and $p\in (0,\IF)$ locally uniformly for $t>0$
\BQNY \pb{\X_I\le   \AIJ+ (p a_n)^{1/p}\x_I \bigl \lvert \absXD = 1- a_n t }  &\to &
\pb{\AII  \bag \VVI  -  1_p \AIJ  \le  t^{-1/p}\x_I },  \quad n\to \IF
\EQNY
holds. Hence for any  $x_d\in (-\IF,0)$
\BQNY
\lefteqn{\pb{ \X_I    \le  \AIJ+ (p a_n)^{1/p} \x_I, (\absXD-1)/a_n \le x_d \bigl \lvert \absXD > 1- a_n  }}\\
 &=&
\int_{-1}^{x_d}\pb{\X_I\le  \AIJ  +  (p a_n)^{1/p} \x_I \bigl \lvert \absXD
= 1+ a_n t }\, d H(1+ a_n t )/
\overline H(1- a_n)  \\
 &\to & \int_{-1}^{x_d} \pk{\AII \bag  \VVI - 1_p \AIJ \le \abs{t}^{-1/p}\x_I}\, d ( 1- \abs{t}^{\alpha+\gamma} ), \quad \ntoi \\
&=& \pk{\abs{\Yy}^{1/p}(\AII \bag \VVI - 1_p \AIJ  )\le \x_I, \Yy \le x_d},
\EQNY
with $\Yy$ a negative random variable with distribution function $M(s):= 1- \abs{s}^{\alpha+\gamma}, s\in [-1,0]$.
Thus the result follows.

c)  Set $\vk{u}_{n}:=  \vk{a}/a_n, n\ge 1,c:=\normAAJ$ and define a positive random variable $R_n,n\ge 1$ as in \eqref{eq:ER}.
Since $W_p$ is bounded for any $\lambda\in (0,\IF)$ we have $\E{W_p^{\lambda}}\in (0,\IF)$. Hence  for any $z>0$ \nelem{lem:1} implies
\BQNY \limit{n}\pk{ R_n > a_n z } &=&  \limit{n}\frac{
\int_{a_n(c^p + z^p)^{1/p}}^\IF g(a_n c/r)\frac{1}{r} \, d F(r)}
 {\int_{a_n c }^\IF g(a_n c /r)\frac{1}{r} \, d F(r)}\\
 &=&  \frac{\int_{ (c^p + z^p)^{1/p}}^\IF g(c /r) r^{-\gamma-2}\, dr}
 { \int_{ c}^\IF g(c /r) r^{-\gamma-2}\, dr}\\
&=& \pk{\kal{R}> z},
\EQNY
where $\kal{R}\sim Q_{M,g,c}$ being independent of  $\VVI$
and $M(x):= 1- c^{\gamma} x^{-\gamma}, \forall x\ge c>0.$ Thus the result follows easily. \QED

\COM{ \BQNY \limit{n}\pb{ \frac{1}{a_n} R_{a_n, \alpha, \beta,
\delta}> z } &=&  \limit{n}\frac{  \int_{ a_n(1 + z^p)^{1/p}}^\IF
[s^{\pst}- a_n^{\pst}]^{\alpha-1} s^{-r}\, d F(s)}
 {\int_{a_n}^\IF
[s^{\pst}- a_n ^{\pst}]^{\alpha-1} s^{-r}\, d F(s)}\\
&=&  \frac{ \int_{ (1 + z^p)^{1/p}}^\infty[s^{\pst}- 1]^{\alpha-1}
s^{-r- \gamma -1}\, d s} { \int_{ 1 }^\infty
[s^{\pst}- 1]^{\alpha-1} s^{-r- \gamma -1}\, d s}.
\EQNY
Consequently as $n\to \infty$ we obtain
\BQNY
\frac{1}{\norm{\uunJ}_{A,p}} \AII R_{\norm{\uunJ}_{A,p}, \alpha, \beta,
\delta}  \VVI &\todis & \AII \kal{R}_{\alpha,\beta,\gamma,p,\delta} \VVI ,
\EQNY
where $\kal{R}_{\alpha,\beta,\gamma,p,\delta}>0$ almost surely being further independent of $\VVI$
with distribution function
$$
 \frac{ \int_{ 1}^{(1 + x^p)^{1/p}}[s^{\pst}- 1]^{\alpha-1} s^{-r- \gamma -1}\, d s}
{ \int_{ 1 }^\infty [s^{\pst}- 1]^{\alpha-1} s^{-r- \gamma -1}\, d s}, \quad x>0, $$
 thus the proof is complete. }

\prooftheo{theo:coco} The proof can be established by extending Result 1 of Nagaraja and David (1994) to the higher dimensional setup and utilising further the Kotz approximation.\\
We give next the sketch of another proof. For notational  simplicity assume that $k=2$. Define next
$b_n:= H^{-1}(1- 1/n), a_n:=1/w(b_n),n\ge 1$ and let $A_n,B_n$  be as in \eqref{ANBN}.
 In view of \eqref{eq:concos} for $n>1$ we have
\BQNY
\pb{\frac{Y_{[n:n]}- B_n}{A_n}\le y_1, \frac{Y_{[n-1:n]}- B_n}{A_n}  \le y_2}
&=& \int_{ x_1> x_2} \prod_{i=1}^2\pk{Y_1 \le A_ny_i+B_n \lvert X_1 = a_nx_i+b_n}\,  d D_n^*(x_1,x_2),
\EQNY
with $D_n^*(x_1,x_2)= \pk{X_{n:n} \le a_n x_1+b_n, X_{n-1:n} \le a_n x_2 + b_n}$.
As in the proof of Theorem 3.1 we have that the distribution function $H$ of $X_1$ is in the Gumbel max-domain of attraction with scaling function $w$, implying thus the joint convergence of upper order statistics (see e.g., Falk et al.\ (2004)), i.e.,
$$ \limit{n} D_n^*(x_1,x_2) = \kal{D}(x_1,x_2), \quad \forall x_1,x_2, \quad x_1>x_2,$$
with $\kal{D}$ a bivariate distribution function. By the assumptions and the properties of the scaling function $w$
$$ \prod_{i=1}^2\pk{Y_1 \le A_ny_i+B_n \lvert X_1=a_nx_i+b_n} \to   \prod_{i=1}^2\pk{ I_2 \kal{R}_\alpha \le y_i}, \quad n\to \IF$$
holds locally uniformly for $x_1,x_2 \inr$. Hence with similar arguments as in the proof of Theorem 4.1 in Hashorva (2008b) we obtain
\BQNY
\limit{n}\pb{\frac{Y_{[n:n]}- B_n}{A_n}\le y_1, \frac{Y_{[n-1:n]}- B_n}{A_n}  \le y_2}&=&
\prod_{i=1}^2\pk{ I_2 \kal{R}_\alpha \le y_i}
 \int_{ x_1> x_2} \,  d \kal{D}(x_1,x_2)=\prod_{i=1}^2\pk{ I_2 \kal{R}_\alpha \le y_i}  ,
\EQNY
thus the result follows.
\QED

{\bf Acknowledgement:} I would like to thank Professor Samuel Kotz for inspiring discussions, and the Referee for
detailed review, several corrections and valuable suggestions.

\bibliographystyle{plain}

\end{document}